%% file: RR-6788.tex
\documentclass[a4paper]{article}
\usepackage{RR}
\usepackage{hyperref}
\RRdate{D\'ecembre 2008}
\usepackage{amsmath}
\usepackage{amssymb}
\usepackage{graphicx}
\usepackage{epsfig}
\usepackage{color} 

\newcommand{\CC}{ \mathcal{C}}

\newcommand{\cqfd}{\mbox{}\hfill\rule{.8em}{1.8ex}}

\newcommand{\egaldef}{\stackrel{\mathrm{def}}{=}}

\newcommand{\rf}{{\rm ref}}

\newcommand{\mn}{{\rm min}}
\newcommand{\mx}{{\rm max}}

\newcommand{\ttt}{{\mathbb{T}}}

\newtheorem{Remark}{Remark}

\def\bbbr{{I\!\!R}}

\RRauthor{
Guy Chavent
\thanks{INRIA Paris-Rocquencourt, and Ceremade, Paris-Dauphine university, France,
E-mail: guy.chavent@inria.fr }}
\authorhead{G. Chavent}
\RRtitle{Une formulation des \'ecoulements triphasiques en pression globale et totalement \'equivalente}
\RRetitle{A Fully Equivalent Global Pressure Formulation for Three-Phase Compressible Flow}
\titlehead{A Fully Equivalent Global Pressure Formulation}

\RRabstract{
We introduce a new global pressure formulation for immiscible three-phase compressible flows in porous media which is fully equivalent to the original equations, unlike the one introduced in \cite{CJ86}. In this formulation, the total volumetric flow of the three fluids and the global pressure follow a classical Darcy law, which simplifies the resolution of the pressure equation.  However, this global pressure formulation exists only for Total Differential (TD) three-phase data, which depend only on two functions of saturations and global pressure: the global capillary pressure and the global mobility.
Hence we introduce a class of interpolation which constructs such TD-three-phase data from any set of three two-phase data (for each pair of fluids) which satisfy a TD-compatibility condition. }

\RRresume{
Nous introduisons une nouvelle formulation en pression globale pour les \'ecoulements triphasiques immiscibles compressibles en milieu poreux. A la diff\'erence de \cite{CJ86}, cette formulation est totalement \'equivalente aux \'equations physiques originales, mais elle permet toujours d'exprimer le flux volum\'etrique total des trois fluides \`a l'aide de la pression globale par une loi de Darcy classique, ce qui simplifie la r\'esolution de l'\'equation en pression. Cette formulation n'existe cependant que pour des perm\'eabilit\'es triphasiques satisfaisant une condition de Diff\'erentielle Totale, qui ne d\'ependent plus que de deux fonctions des saturations et de la pression globale~: la pression capillaire globale  et la mobilit\'e globale. C'est pourquoi nous introduisons une famille d'interpolation permettant de construire de telles perm\'eabilit\'es \`a partir de la connaissance de trois  jeu de donn\'ees diphasiques (un pour chaque paire de fluides) v\'erifiant une condition de compatibilit\'e.}

\RRmotcle{Ecoulement triphasique en milieu poreux.}
\RRkeyword{Three-phase flow in porous media, global pressure.}

\RRprojets{Estime}
\RRtheme{\THNum}
\RCParis

\begin{document}
\RRNo{6788}
\makeRR

\section{Introduction}
\label{sect intro}
\setcounter{equation}{0}

The numerical simulation of three-phase immiscible compressible flows in porous media requires the knowledge of three-phase relative permeabilities and capillary pressures. In practice, experimental values for these data are available only on the boundary of the ternary diagram, i. e. for the three two phase flows corresponding to each pair of fluids. 
The three-phase relative permeabilities are then derived by interpolation from these sets of two-phase data (see for example \cite{stone}). The existence of various interpolation formula is in itself a sign that
none of these formula is the ultimate truth.

We take advantage in this paper of this lack of experimental three-phase data, and introduce a new \emph{class of TD-interpolations},  which is designed to simplify the numerical simulation of the flow by allowing the use of a \emph{global pressure formulation}. 

The global pressure reformulation of the original flow equation was introduced for incompressible two-phase flows in \cite{C76,AM78}, and generalized, under 
the approximation that the volume factors be evaluated at the new global pressure instead of the corresponding phase pressure, to \emph{compressible} two and three-phase flows in \cite{CJ86}. For three-phase flows (both incompressible  and compressible), the original global pressure formulation required that the thee-phase relative permeabilities and capillary pressures satisfy a \emph{Total Differential (TD)-condition}. An algorithm for the determination of such \emph{TD-three-phase data} was given in  \cite{CJ86} and implemented in \cite{CS85,jegou97}, but its complexity has limited the use of the global pressure in numerical simulation codes. Nevertheless, comparison with other approaches \cite{CE97} show the computational effectiveness of this approach when it can be put to work, which may explain the current revival of interest for the global pressure~: the approximation on the volume factors required in the original global pressure formulation of compressible flows 
has been lifted recently for two phase flows in \cite{AJ08} in an independantly work, and is lifted for  three-phase flows in the present paper.

So we first show that, under an integral \emph{Total Differential (TD)-condition} linking three-phase relative permeabilities and capillary pressures, the total volumetric flow of the three phases is governed by a single phase Darcy-like law for a \emph{new global pressure} variable $P$. This new formulation is fully equivalent to the original three-phase compressible equations, in opposition to the formulation proposed in \cite{CJ86}. Three-phase relative permeabilities and capillary pressures which satisfy the TD-condition we be called \emph{TD-three-phase data}.

Next we show that any collection of three two-phase data sets (defined on the boundary $\partial \ttt$ of the ternary diagram $\ttt$) which is the trace of some TD-three-phase data (defined on $\ttt$) satisfies necessarily a \emph{TD-compatibility condition}. 

Finally, we study the \emph{TD-interpolation} of a collection of three two-phase data sets which satisfies the TD-compatibility condition~:   we show that TD-interpolation reduces to the choice of two functions over the ternary diagram (one global capillary pressure, to be defined in section \ref{sect searching for p}, and one global mobility) which satisfy boundary conditions determined by the three given two-phase data sets. This is expected to make more easy the numerical  determination of TD-three-phase data.

\section{The three-phase immiscible compressible equations }
\label{sect three phase equations}
\setcounter{equation}{0}

Let the fluids be numbered in order of decreasing wettability, and denote by upper case letters $S_{j}, P_{j}, \dots$ saturation and pressure distributions (function of the space and time variables $x,t$), and by lower case letters $s_{j},p_{j}, \dots$ saturation and pressure levels (real positive numbers). It will be convenient to use vector notations for the saturations~:
\begin{equation}
\label{eq 0}
\left\{\begin{array}{lll}
        	S=(S_{1},S_{3})&=&vector \  of \  saturation \  fields \ , \\
	s=(s_{1},s_{3})&=&vector \  of \  saturation \  levels \ , 
  \end{array}\right.
\end{equation}
and to denote by $\ttt$ and $\partial \ttt$ the ternary diagram and its boundary~:
\begin{equation}
\label{eq 0bis}
\left\{\begin{array}{lll}
    \ttt & = & \{ s \mbox{ such that }    0\leq s_{j}\leq 1\ , \ j=1,3\ ,\ s_{1}+s_{3}\leq1 \}  \\
  \partial \ttt & = & \{Ês \in \ttt  \mbox{ such that } s_{1}=0 \mbox{ or } s_{3}=0 \mbox{ or } s_{1}
+s_{3}=1 \}
\end{array}\right.
\end{equation}

\subsection{Conservation laws~:}
  For each phase $j=1$ (water), $2$ (oil), $3$ (gas) one has~: 
  \begin{equation}
\label{eq 1}
  \frac{\partial}{\partial t} \big\{ \phi(x,P_{\mathrm{pore}}) \, B_{j}(P_{j})\, S_{j} \big\}+\nabla \cdot \varphi_{j}=0 \quad,\quad j=1,2,3.
\end{equation}
where:
\begin{equation}
\label{eq 1bis}  \left\{ \begin{array}{lll}
       x    & =& vector \ of \ spatial \ coordinates, \\
       \phi(x,p) &=& porosity \  at \  location \  x \ and \ pressure \  p,    \\
         P_{pore}&=& pore \  pressure,  \ \simeq P_{1},P_{2},\ or \ P_{3}.
\end{array} \right.
\end{equation}
and where, for each phase $j=1,2,3$ :
\begin{equation}
\label{eq 1ter}  \left\{ \begin{array}{lllll }
       P_{j}  &=&P_{j}(x,t)&=  pressure, \\
       S_{j}&=&S_{j}(x,t)&= reduced \ saturation\ , \ S=S_{1},S_{2},S_{3} \\
       &&& \hspace{1,3em} 0\leq S_{j}\leq 1\ ,\ S_{1}+S_{2}+S_{3}=1\ , \hspace{-2em}\\
       \varphi_{j}&=&\varphi_{j}(x,t)&=volumetric \ flow \ vector \ at \ reference \ pressure .\hspace{-2em} \\
       B_{j}(p_{j})&= &\rho_{j}(p_{j})/ \rho_{j}^\rf &=volume\  factor \ at \ pressure \ p_{j},   \\
       where: &&\rho_{j}(p_{j})&=density, \ (mass \ per \ unit \ volume)  \ at \ pressure \ p_{j} \hspace{-2em}\\
       && \rho_{j}^\rf &=density \ at \ reference \ pressure, \hspace{-20em} 
       \end{array} \right.
\end{equation}
  \subsection{Muskat law~:} The volumetric flow vector of each phase  $j=1,2,3$ at reference pressure  is given by~:
  \begin{equation}
\label{eq 2}
\varphi_{j}=-K(x) \,d_{j}(P_{j})\,kr_{j}(S)(\nabla P_{j}-\rho_{j}(P_{j})g \nabla Z)
\end{equation}
where:
\begin{equation}
\label{eq 2bis}  \left\{ \begin{array}{lll}
        K(x)&= &absolute \ permeability  \ at \ location \ x,    \\
        d_{j}(p_{j})& =& B_{j}/\mu_{j}=phase \ mobility  \ at \ pressure \ p_{j}, \\
        \mu_{j}(p_{j})&= &phase \ viscosity  \ at \ pressure \ p_{j}, \\
        kr_{j}(s)& =& phase \ relative \ permeability \ at \ saturation \ levels \  s, \hspace{-2em}\\
                g &=& gravity \ constant \ , \\
        Z(x)& =&depth \ of \ location \ x.
\end{array} \right.
\end{equation}
  \subsection{Capillary pressure law~:}
\begin{equation}
\label{eq 3}  \left\{ \begin{array}{lll}
       P_{1}-P_{2} &=& P_{c}^{12}(S_{1})  \ , \\
       P_{3}-P_{2} &=& P_{c}^{32}(S_{3}) \ ,
       \end{array} \right.
\end{equation}
where:
\begin{equation}
\label{eq 3bis}  \left\{ \begin{array}{lll}
       P_{c}^{12}(s_{1})&= &water\!-\!oil \ capillary \ pressure \ at \ water \ saturation \ level \  s_{1},  \hspace{-2em}   \\
       P_{c}^{32}(s_{3})& = &gas\!-\! oil \ capillary \ pressure \ at \ gas \ saturation \ level \  s_{3}. 
       \end{array} \right.
\end{equation}

\section{Classical resolution~: ``pressure equation''}
\label{sect classical resolution pressure equation}
\setcounter{equation}{0}
The numerical resolution of equations (\ref{eq 1})(\ref{eq 2})(\ref{eq 3}) is usually done by solving one ``pressure equation'', to be defined in the next section, with respect of one of the phase pressure, say $P_{2}$, and two of the  ``saturation equations'' (\ref{eq 1}), say with respect to $S_{1}$ and $S_{3}$. 
\subsection{Forming the ``pressure''  equation~:} Summing up the three conservation laws (\ref{eq 1}) gives~:
\begin{equation}
\label{eq 4}
  \frac{\partial}{\partial t} \big\{ \phi(x,P_{\mathrm{pore}}) \sum_{j=1}^3\, B_{j}(P_{j}) \, S_{j} \big\}+\nabla \cdot q=0 \ ,
\end{equation}
where $q$ is the global  volumetric  flow  vector:
\begin{equation}
\label{eq 5}  
   q\egaldef \sum_{j=1}^3 \varphi_{j}= 
   	-K\lambda \big\{ \nabla P_{2}+f_{1} \nabla P_{c}^{12}+f_{3} \nabla P_{c}^{13}-\rho g \nabla Z \big\} \ .
\end{equation}
Here $\lambda,f_{1},f_{3},\rho$ are the global mobility, the water and oil fractional flows and the global density expressed as function of the \emph{oil pressure level} $p_{2}$, using the capillary pressure laws (\ref{eq 3})~:
\begin{eqnarray}
     \lambda(s,p_{2}) &=& kr_{1}(s)\,d_{1}(p_{2}+P_{c}^{12}(s_{1}))  +kr_{2}(s)\,d_{2}(p_{2}) +kr_{3}(s)\,d_{3}(p_{2}+P_{c}^{32}(s_{3})) \ , \label{eq 5-1}\\
     f_{j}(s,p_{2}) & = & kr_{j}(s)\,d_{j}(p_{2}+P_{c}^{j2}(s_{j}))/\lambda(s,p_{2}) \ , \  j=1,3 \ ,\label{eq 5-2}\\
     f_{2}(s,p_{2}) & = & 1- f_{1}(s,p_{2}) - f_{3}(s,p_{2}) \ ,\label{eq 5-3}\\
    \rho(s,p_{2})   & = & \sum_{j=1,2,3}  f_{j}(s,p_{2}) \rho_{j}(p_{j}) \ .\label{eq 5-4}
 \end{eqnarray}
In this approach, equations (\ref{eq 4})(\ref{eq 5}) have to be solved for the oil pressure $P_{2}$ at each time step..

\subsection{Difficulties~:} 
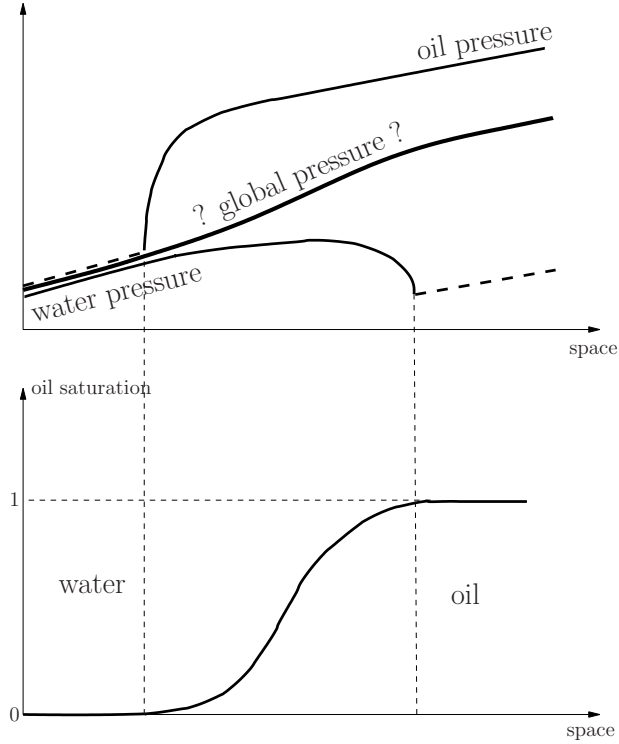
\begin{figure}[ht]
\begin{center}
\begin{minipage}{12.cm}
\centerline{\resizebox{80mm}{!}{\input{figure_1b_B_W.pstex_t}}}
\caption{Behaviour of the individual phase pressure across a front}
\label{fig 1}
\end{minipage}
\end{center}
\end{figure}

Equation (\ref{eq 4})(\ref{eq 5})  is not a classical pressure equation, because of the gradient of capillary pressure terms in (\ref{eq 5}), which makes its discretization delicate. Moreover,
 the individual phase pressures can be singular, as it is shown on figure \ref{fig 1}, for a two-phase problem~: the difference between water and oil pressure is the capillary pressure, whose derivative is infinite at residual saturations - hence the water and oil pressure cannot be both regular across a water-oil front. So the choice of $P_{2}$ as numerical unknown is likely to require a fine or adaptive mesh.

\subsection{Let us have a dream \dots} Could it be possible to find a ``global pressure'' field P such that~:
\begin{itemize}
  \item  $P$ is smooth, in opposition to the individual phase pressures,
  \item $P$ governs the global volumetric flow of the three fluids~:
\begin{equation}
\label{eq 6}
  q=-Kd\big\{ \nabla P-\rho g \nabla Z \big\} \ ,
\end{equation}
where now $d$ and $\rho$ denote the global mobility and density at  a global pressure 
level $p$ - compare with equation (\ref{eq 5-1})(\ref{eq 5-4})~:
 \begin{equation}
\label{eq 6bis}
d(s, p)\egaldef\lambda(s,p_{2})  \quad , \quad \rho(s,p) \egaldef  \rho(s,p_{2})\ ,
\end{equation}
\end{itemize}
(we have kept for simplicity the same notation for the global density as function of $p_{2}$ and~$p$)

\section{Searching for a global pressure P}
\label{sect searching for p}
\setcounter{equation}{0}

Comparison of (\ref{eq 5})(\ref{eq 6}) suggests to search for a new dependant \emph{global pressure} variable $P$ by setting~:
\begin{equation}
\label{eq 7}
  P=P_{2}+P_{cg}(S,P)
\end{equation}
where the \emph{global capillary pressure} function~: 
\begin{equation}
\label{eq 7bis}
   (s,p) \leadsto P_{cg}(s,p)
\end{equation}
is to be chosen such that~:
\begin{itemize}
 \item \textbf{\dots $P$ is in the range of the phase pressures.} With the chosen wettability conventions, the capillary pressures  (\ref{eq 3bis}) satisfy~:
 \begin{equation}
\label{eq 13} \left\{ \begin{array}{ccccc}
P_{c}^{12}(s_{1}) \leq 0 &,&P_{c}^{12}(1) = 0 &,& \displaystyle  \frac{dP_{c}^{12}}{ds_{1}}(s_{1}) \geq 0 \\ 
\\
P_{c}^{32}(s_{3}) \geq 0 &,&P_{c}^{32}(0) = 0 &,& \displaystyle  \frac{dP_{c}^{32}}{ds_{3}}(s_{3}) \geq 0 
\end{array} \right.
\end{equation}
  so that~:
\begin{equation}
\label{eq 13bis}
  P_{\mn} \leq P_{1} \leq P_{2} \leq P_{3} \leq P_{\mx} \quad \mbox{at all points $x$ and time $t$}\ ,
\end{equation}
where $P_{\mn}$ is a known lower bound to the water pressure $P_{1}$, and  $P_{\mx}$ is a known upper bound to the gas pressure $P_{3}$.

One wants that the global pressure satisfies a similar inequality~:
\begin{equation}
\label{eq 14}
  P_{\mn} \leq P_{1} \leq P \leq P_{3} \leq P_{\mx} \quad \mbox{at all points $x$ and time $t$}\ .
\end{equation}
This will make $P$ a natural candidate for the pore pressure $P_{pore}$  required to evaluate the porosity $\phi$ (see (\ref{eq 1bis})).
  \item \textbf{\dots our dream comes (almost \dots) true~!}
Comparison of equations (\ref{eq 5}) and (\ref{eq 6}) suggests to require that~:
  \begin{equation}
\label{eq 10} \hspace*{-0,7em}\left\{ \begin{array}{l}
  \mbox{for any saturation/global pressure fields $S_{1}(x,t),S_{3}(x,t),P(x,t)$, }       \\
  \mbox{the global capillary pressure function (\ref{eq 7bis}) satisfies:}   \\
  \nabla P_{cg}(S,P) = f_{1}(S,\underbrace{P-P_{cg}(S,P)}_{\displaystyle P_{2}}) \, \nabla P_{c}^{12}(S_{1}) 
  				+f_{3}(S,\underbrace{P-P_{cg}(S,P)}_{\displaystyle P_{2}}) \, \nabla P_{c}^{32}(S_{3}) \hspace{-2em}\\
   \hspace{15em}+\displaystyle \frac{\partial P_{cg}}{\partial p}(S,P)\, \nabla P \ ,
\end{array} \right.
\end{equation}
where $\nabla$ denotes, as in equations (\ref{eq 2})(\ref{eq 5}) (\ref{eq 6}), the gradient with respect to the space variables $x$. When this condition is satisfied, equation (\ref{eq 5}) becomes, with the notations (\ref{eq 6bis})~:
\begin{equation}
\label{eq 10bis}
    q=-Kd\big\{ (1-\frac{\partial P_{cg}}{\partial p})\nabla P-\rho g \nabla Z \big\} \ ,
\end{equation}
which is the desired equation (\ref{eq 6}) up to the factor $(1-\partial P_{cg}/\partial P)$~!
  \item \textbf{\dots $P$  is uniquely determined by (\ref{eq 7}) once $P_{2}$ and $S=(S_{1},S_{3})$ are known.}  This property will make it possible to use $(P,S)$ as dependant variables instead of $(P_{2},S)$. It will be satisfied as soon as the \emph{stability condition}~:
  \begin{equation}
\label{eq 8}
  | \frac{\partial P_{cg}}{\partial p}(s,p) | <1 \ \mbox{ for all } \  s \in \ttt\ ,\ P_\mn \leq p \leq P_{\mx} 
\end{equation}
is satisfied. Notice that the computation of $P$ from $S,P_{2}$ is required only for the determination of the initial global pressure distribution; after that, during the course of the computation, the numerical code knows only $P$, and it is $P_{2}$ which is computed from $P$ and $S$ using (\ref{eq 7}).
\end{itemize}
Hence a global pressure formulation will exist as soon as one can find a \emph{global capillary pressure function} $s,p \leadsto P_{cg}(s,p)$ which satisfies the three conditions (\ref{eq 14})(\ref{eq 10})(\ref{eq 8}).
\subsection{Satisfying condition (\ref{eq 14})~:} Condition (\ref{eq 10}) will constrain the global capillary pressure up to a constant, so we can fix its value at one point of the ternary diagram. Choosing~:
\begin{equation}
\label{eq 15}
  P_{cg}(1,0, p) = 0 \quad \mbox{ for} \ P_\mn \leq p \leq P_{\mx}\ ,
\end{equation}
will imply the desired property (\ref{eq 14}) as it follows immediately from (\ref{eq 16}) below.

\subsection{Satisfying condition (\ref{eq 10})~:}
   Taking successively $S_{1}(x,t)=\mbox{constant}$ and $S_{3}(x,t)=\mbox{constant}$ shows that (\ref{eq 10}) will be satisfied if and only if the global capillary function satisfies~:
\begin{equation}
\label{eq 11} \left\{ \begin{array}{lcl}
       \displaystyle \frac{\partial P_{cg}}{\partial s_{1}}(s,p) & =  &
     	f_{1}(s,p-P_{cg}(s,p)) \, \displaystyle \frac{dP_{c}^{12}}{ds_{1}}(s_{1}) \ ,\\
      \displaystyle \frac{\partial P_{cg}}{\partial s_{3}}(s,p) &  = &
     	f_{3}(s,p-P_{cg}(s,p))  \, \displaystyle \frac{dP_{c}^{32}}{ds_{3}}(s_{3}) \ ,
\end{array} \right.
\end{equation}
for all $s \in \ttt\ ,\ P_\mn \leq p \leq P_{\mx}$. 

This relation implies that, if it exists, $P_{cg}$ satisfies a differential equation along any smooth curve 
$\CC$ of $\ttt$. Hence the existence of a global capillary pressure function $P_{cg}$ satisfying (\ref{eq 11}) amounts to the following \emph{Total Differential (TD)} condition~:
\begin{equation}
\label{eq 11-1}
\left\{\begin{array}{l}
\mbox{for any $s\in \ttt$, $P_\mn \leq p \leq P_{\mx}$ and any smooth curve $\CC:[0,1] \leadsto \ttt$, such} \\
\mbox{that $\CC(0)=(1,0)\ ,\ \CC(1)=s$, the solution of the differential equation~: }\\
   \displaystyle  \frac{d\beta}{dt}=f_{1}(\CC,p-\beta) \, \displaystyle \frac{dP_{c}^{12}}{ds_{1}}(\CC_{1})\,\CC^\prime_{1}
    +  f_{3}(\CC,p-\beta) \, \frac{dP_{c}^{32}}{ds_{3}}(\CC_{3})\,\CC^\prime_{3} \ ,\ \beta(0) = 0 \\
    \mbox{satisfies $\beta(1)=s$ independantly of $\CC$.}
\end{array} \right.
\end{equation}

When this TD-condition is satisfied, the global capillary pressure satisfying (\ref{eq 14})(\ref{eq 10}) is given by~:
\begin{equation}
\label{eq 11-2}
P_{cg}(s,p)=\beta(1).
\end{equation}
\begin{figure}[ht]
\centerline{
\resizebox{70mm}{!}{\input{figure_3_B_W.pstex_t}}
}
\caption{The Total Differential condition for three phase data~: $P_{cg}(s,p)$ is required to be independant of the curve $\CC$ along which it is computed.}
\label{fig 2}
\end{figure}
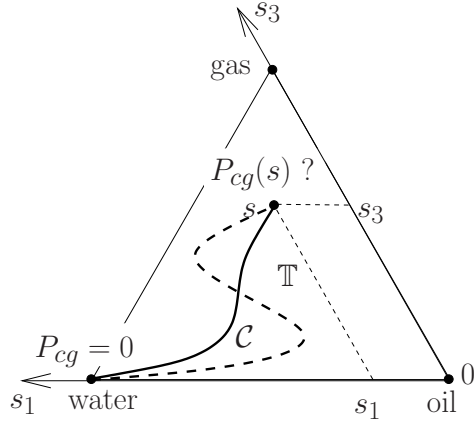

 An \emph{alternative formulation of the TD-condition} is obtained by requiring that  the second derivative 
 $\partial^2 P_{cg}/\partial S_{1}\partial S_{3}$ takes the same value when computed from either line of (\ref{eq 11})~:
\begin{equation}
\label{eq 12}
\left\{\begin{array}{l}
     \mbox{There exists a function $s,p \leadsto P_{cg}(s,p)$ s.t.  for all $s \in \ttt\ ,\ P_\mn \leq p \leq P_{\mx}$ ~:}     \\
   \displaystyle \frac{\partial}{\partial s_{3}}\big( f_{1}(s,p-P_{cg}(s,p)) \big)
  \displaystyle  \frac{dP_{c}^{12}}{ds_{1}}(s_{1}) =
   \displaystyle \frac{\partial}{\partial s_{1}} \big(f_{3}(s,p-P_{cg}(s,p)) \big)
    \displaystyle  \frac{dP_{c}^{32}}{ds_{3}}(s_{3}) \ .
  \end{array} \right.
\end{equation}
In the approximate global pressure formulation (see \cite{CJ86}), the $P_{cg}$ terms were neglected in the right-hand sides of (\ref{eq 11}), and hence in (\ref{eq 12}). So condition (\ref{eq 12}) could be used to check wether a given set of three-phase data satisfied the TD-condition or not, and to 
 construct fractional flow functions satisfying this condition \cite{CS85}. This approach cannot be used anymore in the exact formulation given here, where the condition  (\ref{eq 12}) begins with ``there exists $P_{cg}$ such that \dots''~:
 one has to resort to condition (\ref{eq 11-1}) to check if a given set of three-phase data satisfy the TD-condition, and to use a more direct - and hopefully simpler - approach for the construction of TD-three-phase data, which we decribe now.

The TD-condition (\ref{eq 11-1}) or (\ref{eq 12}) has the effect of reducing the number of functions which can be chosen freely over the ternary diagram: instead three relative permeabilities $kr_{1},kr_{2},kr_{3}$ which are function of  saturations (or equivalently two fractional flows  $f_{1},f_{3}$ and the global mobility $\lambda$  which are function of  saturations $s=(s_{1},s_{3})$ \emph{and} oil pressure  $p_{2}$), one can choose only in a global pressure formulation
the \emph{global capillary pressure} $P_{cg}$ and the \emph{global mobility} $d$  which are function of saturation $s=(s_{1},s_{3})$ \emph{and} global pressure~$p$.
%

Based on (\ref{eq 11}) the fractional flows and relative permeabilities associated to $P_{cg}$ and $d$ are~: 
\begin{eqnarray}
   \nu_{j}(s,{p}) &=&  {\partial {P_{cg}}}/{\partial s_{j}}(s,{p})
     			\big{/} \ {dP_{c}^{j2}}/{ds_{j}}(s_{j}) \quad, \quad j=1,3 \label{eq 12-1} \\
    \nu_{2}(s,{p}) &=& 1- \nu_{1}(s,{p})  - \nu_{3}(s,{p})  \label{eq 12-2} \\
          kr_{j}(s,{p}) &=& {\nu_{j}(s,{p})}{d}(s,{p})\big/ 
              {d_{j}\big({p}-{P_{cg}}(s,{p})+P_{c}^{j2}(s_{j})\big)} \quad j=1,3 \ ,    \label{eq 12-3}\\
     kr_{2}(s,{p}) &=&  {(1-\nu_{1}(s,{p})-\nu_{3}(s,{p}))} {d}(s,{p})\big/ 
             {d_{2}\big(p-{P_{cg}}(s,{p})\big)} \ ,    \label{eq 12-4}
\end{eqnarray}
where we have used the notation $\nu_{j}$ for the fractional flows considered as \emph{function of the
 global pressure} $p$. They are related to the fractional flows $f_{j}$ as function of the oil pressure $p_{2}$ by:
\begin{equation}
\label{eq 12-5}
    \nu_{j}(s,p)= f_{j}(s,p-P_{cg}(s,p)) \quad , \quad j=1,2,3 \ .
\end{equation}

\subsection{Satisfying condition (\ref{eq 8})~:}
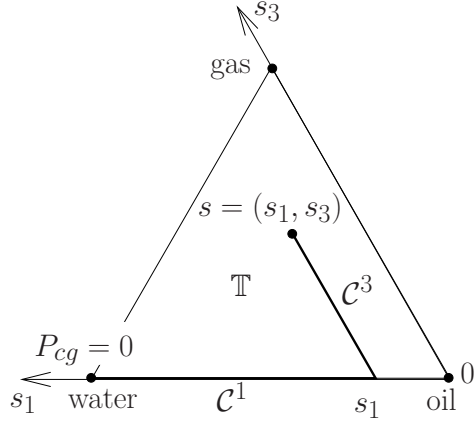
\begin{figure}[ht]
\centerline{
\resizebox{70mm}{!}{\input{figure_3a_B_W.pstex_t}}
}
\caption{One convenient path for the computation of  $\partial P_{cg}/\partial p$.}
\label{fig 3}
\end{figure}

Let now the three-phase data satisfy the TD-condition  (\ref{eq 11-1}),  $P_{cg}(s,t)$ be the global capillary pressure defined by (\ref{eq 15})(\ref{eq 11-1})(\ref{eq 11-2}),  $s=(s_{1},s_{3})\in \ttt$, $P_\mn \leq p \leq P_{\mx}$ be given, and $\CC^{1}: t \leadsto (1-t+ts_{1},0)$ and $\CC^{3}: t \leadsto (s_{1}, ts_{3})$ be the curves shown in figure \ref{fig 3}. 
The functions $\beta^{j}(t)=P_{cg}(\CC^j(t),p)\ ,\  j= 1,3$ satisfy the differential equations (see (\ref{eq 11-1}))~:
\begin{equation}
\label{eq 16}
\left\{\begin{array}{l}
   \displaystyle  \frac{d\beta^1}{dt}= - f_{1}(1-t+ts_{1},0,p-\beta^1) \, \displaystyle \frac{dP_{c}^{12}}{ds_{1}}(1-t+ts_{1})\,(1-s_{1}) \ ,\ \beta^1(0) = 0  \ ,\\
  \displaystyle  \frac{d\beta^3}{dt}= + f_{3}(s_{1}, ts_{3},p-\beta^3) \, \frac{dP_{c}^{32}}{ds_{3}}(ts_{3})\,s_{3} \ ,\ \beta^3(0) = \beta^1(1) \ .\\
 \end{array} \right.
\end{equation}

Derivation of (\ref{eq 16}) with respect to $p$ shows that the functions $\gamma_{j}(t)=\partial P_{cg}/\partial p(\CC_{j}(t),p)$ are given by the the differential equations~:
\begin{equation}
\label{eq 17}
\hspace*{-0,2em} \left\{\begin{array}{l}
   \displaystyle  \frac{d\gamma^1}{dt}= - \frac{\partial f_{1}}{\partial p_{2}}(1-t+ts_{1},0,p-\beta^1)
       \, \displaystyle \frac{dP_{c}^{12}}{ds_{1}}(1-t+ts_{1})\,(1-s_{1}) (1-\gamma^1) \ ,\ \gamma^1(0) = 0  \ ,
       \hspace*{-2em}\\
  \displaystyle  \frac{d\gamma^3}{dt}= + \frac{\partial f_{3}}{\partial p_{2}}(s_{1}, ts_{3},p-\beta^3)
       \, \frac{dP_{c}^{32}}{ds_{3}}(ts_{3})\,s_{3}  (1-\gamma^3) \ ,\ \gamma^3(0) = \gamma^1(1) \ . \\
 \end{array} \right.
\end{equation}
Both equations in (\ref{eq 17}) are of the form~:
\begin{equation}
\label{eq 18}
    \frac{d\gamma}{dt}= \alpha(t)(1-\gamma)\ ,\ \gamma(0)=\gamma_{0} <1  \quad \Longleftrightarrow \quad 
    \gamma(t) = 1-(1-\gamma_{0})\exp{\!\{-\int_{0}^t \alpha(\tau)d\tau\}} \ ,
\end{equation}
so that~:
\begin{equation}
\label{eq 19}
   \frac{\partial P_{cg}}{\partial p}(s,p) = \gamma^3(1) = 1 - \exp{\! \{-\int_{0}^1( \alpha^1(\tau)+\alpha^3(\tau))d\tau \}} \ ,
\end{equation}
which shows that \emph{the coefficient $1-{\partial P_{cg}}/{\partial p}$ in the global Darcy law (\ref{eq 10bis}) is always strictly positive}.

The \emph{stability condition (\ref{eq 8})} that  $|\partial P_{cg}/\partial p| < 1$ is then equivalent to~:
\begin{equation}
\label{eq 20}
    \forall s \in \ttt \ ,\  \forall p \in [P_{\mn}, P_{\mx}] \quad, \quad  \exp{\! \{-\int_{0}^1( \alpha^1(\tau)+\alpha^3(\tau))d\tau \}} < 2 \ .
\end{equation}
Derivation of the fractional flows $f_{1}(s,p_{2})$ and $f_{3}(s,p_{2})$ given by (\ref{eq 5-1})(\ref{eq 5-2}) with respect to $p_{2}$ gives the following expressions for $\alpha^{1}$ and $\alpha^{3}$ (remember~: $kr_{3}=0$ on $\CC^1$)~:
\begin{eqnarray}
\label{eq 20-1}
   \alpha^{1}(t) &=& - \displaystyle \, \frac{kr_{1}d_{1}kr_{2}d_{2}\big( \displaystyle  \frac{d^\prime_{1}}{d_{1}}
              - \displaystyle   \frac{d^\prime_{2}}{d_{2}}\big)}{(kr_{1}d_{1}+kr_{2}d_{2})^2}
              \times  \frac{dP_{c}^{12}}{ds_{1}} \times (1-s_{1})\ , \\
\label{eq 20-2}
  \alpha^{3}(t) &=& + \displaystyle \, \frac{kr_{3}d_{3}kr_{1}d_{1}\big(\displaystyle  \frac{d^\prime_{3}}{d_{3}}
            -\displaystyle  \frac{d^\prime_{1}}{d_{1}}\big)
                                          + kr_{3}d_{3}kr_{2}d_{2}\big(\displaystyle  \frac{d^\prime_{3}}{d_{3}}
            -\displaystyle  \frac{d^\prime_{2}}{d_{2}} \big)}
      {(kr_{1}d_{1}+kr_{2}d_{2}+kr_{3}d_{3})^2} \times  \frac{dP_{c}^{32}}{ds_{3}} \times s_{3} \ ,
\end{eqnarray}
where $kr_{j},d_{j},d^\prime_{j}\, ,\  j=1,2,3$ and $P_{c}^{j2}\,,\, j=1,3$  are evaluated at saturations and phase pressures associated to the point $\CC^1(t)$ (for $\alpha_{1}$) or $\CC^2(t)$ (for $\alpha_{2}$) 
of $\ttt$ at the given global pressure level $p$.


When the fluid densities satisfy~:
\begin{equation}
\label{eq 21}
 \left\{ \begin{array}{l}
   \forall s \in \ttt \ ,\  \forall p \in [P_{\mn}, P_{\mx}]  õ: \\
   \displaystyle \frac{d^\prime_{3}}{d_{3}}(p-P_{cg}(s,p)+P_{c}^{32}(s_{3})) \geq
 \frac{d^\prime_{2}}{d_{2}}(p-P_{cg}(s,p)) \geq
      \frac{d^\prime_{1}}{d_{1}}(p-P_{cg}(s,p)+P_{c}^{12}(s_{1})) \ , \hspace{-1em}
\end{array} \right.
\end{equation}
 $\alpha^1(\tau)$ and $\alpha^3(\tau)$ are positive, so that~: 
\begin{equation}
\label{eq 22}
   \forall s \in \ttt \ ,\  \forall p \in [P_{\mn}, P_{\mx}]  \quad , \quad 0 \leq \frac{\partial P_{cg}}{\partial p}(s,p) < 1 \ ,
\end{equation}
and the stability condition (\ref{eq 8}) is satisfied.

Naturally, (\ref{eq 21}) is only a sufficient conditions for  (\ref{eq 8}). In practice, implementation of the global pressure formulation requires the determination of  $\partial P_{cg}/\partial p$ over the ternary diagram $\ttt$ and the global pressure range of interest $[P_{\mn},P_{\mx}]$, for example by formula (\ref{eq 19})(\ref{eq 20-1})(\ref{eq 20-2}), so one 
can check directly on $\partial P_{cg}/\partial p$ if condition (\ref{eq 8}) is satisfied. In the case where condition (\ref{eq 8}) is not satisfied, it is not known wether this limitation is due to the global pressure formulation, or is the sign that the original equations are not anymore well-posed, or change  of type.

\begin{Remark}
\label{rem 23}
In the case of a \emph{water-oil-gas system}, one has $d^\prime_{1}=d^\prime_{2}=0$, and condition 
(\ref{eq 21}) is satisfied~:  \emph{the global capillary pressure, when it exists, satisfies always (\ref{eq 22})}. In that case, $\partial P_{cg}/\partial p=0$ on the water - oil side of $\ttt$, so that $\alpha^1=0$ and (\ref{eq 19}) simplifies somewhat.
\cqfd
\end{Remark}

\section{TD-interpolation of two-phase data given on the sides of the ternary diagram $\ttt$ }
\label{sect td-interpolation}
\setcounter{equation}{0}

Let two-phase relative permeabilities and capillary pressure curves be given on the boundary $\partial \ttt$ of  $\ttt$ for each system of  two-fluids~:
\begin{equation}
\label{eq C1}\left\{ \begin{array}{lll}
water-oil~:  &  s_{1} \leadsto kr_{1}^{12}(s_{1}), \  kr_{2}^{12}(s_{1}) & P_{c}^{12}(s_{1}) \ ,  \\
 water-gas~: & s_{1} \leadsto kr_{1}^{13}(s_{1}), \  kr_{3}^{13}(s_{1}) & P_{c}^{13}(s_{1})=  P_{c}^{12}(s_{1})-P_{c}^{32}(1-s_{1}) \ , \\
  gas-oil~:    &  s_{3} \leadsto kr_{3}^{23}(s_{3}), \  kr_{2}^{23}(s_{3}) & P_{c}^{32}(s_{3})  \ , 
\end{array} \right.
\end{equation} 
as well as the three phase mobility functions (see (\ref{eq 2bis}))~:
\begin{equation}
\label{eq C2}
  p_{j} \leadsto d_{j}(p_{j}), \quad j=1,2,3 \ ,
\end{equation}

  The question arises of wether one can find, for each given global pressure level $P_{\mn} \leq p \leq P_{\mx}$, three-phase relative permeabilities $s\in \ttt \leadsto kr_{j}(s,p)\in \bbbr, \ j= 1,2,3$ (compare with (\ref{eq 2bis})) which \dots
  \begin{enumerate}
  \item honor the two-phase data on $\partial \ttt$,
  \item  \emph{and} satisfiy the TD-condition (\ref{eq 11-1}) on $\ttt$.
\end{enumerate}
Satisfying condition 1 and 2 amounts to search for a global capillary pressure function  
$(s,p) \leadsto P_{cg}(s,p)$ and a global mobility function 
$(s,p) \leadsto d(s,p)$
such that  the fractional flow $\nu_{j}(s,p)$ defined by (\ref{eq 11})~:
\begin{equation}
\label{eq C2-2} \left\{ \begin{array}{rcl}
    \nu_{1}(s,p) &=&  \displaystyle \frac{\partial P_{cg}}{\partial s_{1}}(s,p)
     			\Big{/}	\  \displaystyle \frac{dP_{c}^{12}}{ds_{1}}(s_{1}) \ ,\\
    \nu_{3}(s,p) &=&  \displaystyle \frac{\partial P_{cg}}{\partial s_{3}}(s,p)
    		 	\Big{/}	\  \displaystyle \frac{dP_{c}^{32}}{ds_{3}}(s_{3}) \ ,
\end{array} \right.
\end{equation}
produce three-phase relative permeabilities  $kr_{j}(s,p)$~:
\begin{equation}
\label{eq C2-1} \left\{ \begin{array}{lccl}
    kr_{1}(s,p) & = &\displaystyle \frac{\nu_{1}(s,p)}{d_{1}(p-P_{cg}(s,p)+P_{c}^{12}(s_{1}))} & d(s,p)   \ , \\
    kr_{2}(s,p) & = &\displaystyle \frac{1-\nu_{1}(s,p)-\nu_{3}(s,p)}{d_{2}(p-P_{cg}(s,p))} & d(s,p)  \ ,  \\
    kr_{3}(s,p) & = &\displaystyle \frac{\nu_{3}(s,p)}{d_{3}(p-P_{cg}(s,p)+P_{c}^{32}(s_{1}))} &d(s,p) \ ,
\end{array} \right.
\end{equation}
which coincide with the given two phase data (\ref{eq C1}) on $\partial \ttt$.

\subsection{Determination of $P_{cg}(s,p)$ on $\partial \ttt$.}

Let $P_{\mn} \leq p \leq P_{\mx}$ be given, and denote by
\begin{equation}
\label{eq C3-0}
    \CC^{12}(t)=(1-t,0)\quad , \quad \CC^{23}(t)=(0,t) \quad , \quad \CC^{13}(t)=(1-t,t)\quad ,\quad 0 \leq t \leq 1 
\end{equation}
a parameterization of the three edges of $\partial \ttt$. When it exists, the global capillary pressure $\beta^{ij}(t)=P_{cg}(\CC^{ij}(t)), i,j=1,2,3, i<j$ on $\partial \ttt$ satisfies, according to (\ref{eq 11-1})~:
\begin{itemize}
 \item \textbf{along the water-oil-gas sides~:}
\begin{equation}
\label{eq C3-1}
\left\{\begin{array}{l}
   \displaystyle  \frac{d\beta^{12}}{dt}= - \underbrace{f_{1}^{\rm data}(1-t,0,p-\beta^{12})}_{\displaystyle 
                \egaldef \,\nu_{1}^{12,\rm data}(t,p)}  \, \displaystyle \frac{dP_{c}^{12}}{ds_{1}}(1-t) \ ,\ \beta^{12}(0) = 0  \ ,\\
  \displaystyle  \frac{d\beta^{23}}{dt}= + \underbrace{f_{3}^{\rm data}(0, t,p-\beta^{23})}_{\displaystyle \egaldef \,\nu_{3}^{23,\rm data}(t,p)}  \, \frac{dP_{c}^{32}}{ds_{3}}(t)\,s_{3} \ ,\ \beta^{23}(0) = \beta^{12}(1) \ .
 \end{array} \right.
\end{equation}
  \item \textbf{along the water-gas side~:}
 \begin{equation}
\label{eq C3-2}
\left\{\begin{array}{lcl}
    \displaystyle \frac{d\beta^{13}}{dt}&= &- \underbrace{f_{1}^{\rm data}(1-t,t,p-\beta^{13})}_{\displaystyle \egaldef \,\nu_{1}^{13,\rm data}(t,p)}  \, \displaystyle \frac{dP_{c}^{12}}{ds_{1}}(1-t) \\
    &&  + \underbrace{f_{3}^{\rm data}(1-t, t,p-\beta^{13})}_{\displaystyle \egaldef \,\nu_{3}^{13,\rm data}(t,p)} \, \displaystyle \frac{dP_{c}^{32}}{ds_{3}}(t) \ , \\
     \beta^{13}(0) &= &0 \ ,
\end{array}\right.
\end{equation}
 \end{itemize}
 where  $f_{1}^{\rm data}(s,p_{2})$ and $f_{3}^{\rm data}(s,p_{2})$ denote, for $s \in \partial \ttt$, the
two-phase fractional flows as function of the oil pressure derived from the data  (\ref{eq C1})(\ref{eq C2}) using two phase versions of (\ref{eq 5-1})(\ref{eq 5-2}). 


Using the TD-condition (\ref{eq 11-1}), one sees that the existence of a global capillary function implies that the two-phase data (\ref{eq C1})(\ref{eq C2})  satisfy necessarily the \emph{TD-compatibility condition}  
\begin{equation}
\label{eq C3-3}
   \beta^{23}(1) = \beta^{13}(1) \ ,
\end{equation}
i.e., using (\ref{eq C3-1})(\ref{eq C3-2}) and the two-phase fractional flows $\nu_{1}^{12,\rm data},\nu_{3}^{23,\rm data},\nu_{1}^{13,\rm data},\nu_{3}^{13,\rm data}$ defined there~:
\begin{equation}
\label{eq C3-4}
   \int_{0}^1 (\nu_{1}^{12,\rm data}-\nu_{1}^{13,\rm data})  \frac{d P_{c}^{12}}{d s_{1}}
   = \int_{0}^1 (\nu_{3}^{23,\rm data}-\nu_{3}^{13,\rm data})  \frac{d P_{c}^{32}}{d s_{3}} =0
   \  \mbox{ for} \  P_{\mn}\leq p \leq P_{\mx}\ . \hspace{-0,5em}
\end{equation}

This compatibility condition is a constraint on weighted means  of the fractional flows, and hence a non-linear constraint on the mean values of the two-phase relative permeability  data sets  (\ref{eq C1}).
When it is satisfied, one sees that the global capillary function $P_{cg}$ has to satisfy the  \emph{Dirichlet boundary condition}~:
 \begin{equation}
 \label{eq C4}
  P_{cg}=P_{cg}^{\rm data}\  \egaldef \  \left\{ \begin{array}{ll}
      \beta^{12}  &  \mbox{(water-oil edge),}  \\ 
      \beta^{13} & \mbox{(water-gas edge),}  \\
      \beta^{23} &   \mbox{(gas-oil edge),} 
 \end{array} \right.                                                        
\end{equation}

\subsection{Determination of  $\displaystyle  \frac{\partial P_{cg}}{\partial n}(s,p)$ on $\partial \ttt$.}

Now that $P_{cg}$ is known on $\partial \ttt$, equation (\ref{eq 11})  gives, for $s \in \partial \ttt$~:
 \begin{equation}
\label{eq C7} \left\{ \begin{array}{lcl}
   \displaystyle  \frac{\partial P_{cg}}{\partial s_{1}}(s,p) &=&
                f_{1}^{\rm data}(s,p-P_{cg}^{\rm data}(s,p))\displaystyle  \frac{d P_{c}^{12}}{d s_{1}}(s_{1})  \ ,\\
    \displaystyle  \frac{\partial P_{cg}}{\partial s_{3}}(s,p) &=&
                f_{3}^{\rm data}(s,p-P_{cg}^{\rm data}(s,p))\displaystyle  \frac{d P_{c}^{32}}{d s_{3}}(s_{3})  \ .
\end{array} \right.
\end{equation}
Let then $n$ be the unit outer normal to $\ttt$. The normal derivative of $P_{cg}$ on $\partial \ttt$ is given by~:
\begin{equation}
\label{eq C6} \left\{ \begin{array}{lccclcll}
    \displaystyle  \frac{\partial P_{cg}}{\partial n} &=& \frac{\sqrt{3}}{3} (\displaystyle   \frac{\partial P_{cg}}{\partial s_{1}} - 2 \displaystyle  \frac{\partial P_{cg}}{\partial s_{3}}) 
	& \mbox{(water-oil edge)} \ ,\\
     \displaystyle  \frac{\partial P_{cg}}{\partial n} &=& \frac{\sqrt{3}}{3} ( \displaystyle   \frac{\partial P_{cg}}{\partial s_{1}} +
     	\displaystyle  \frac{\partial P_{cg}}{\partial s_{3}})  
	& \mbox{(water-gas edge)} \ , \\
	   \displaystyle  \frac{\partial P_{cg}}{\partial n}&=& \frac{\sqrt{3}}{3}( \displaystyle  \frac{\partial P_{cg}}{\partial s_{3}}  -  2 \displaystyle   \frac{\partial P_{cg}}{\partial s_{1}})
     	 & \mbox{(gas-oil edge)} \ ,

\end{array} \right.
\end{equation}
Combining (\ref{eq C7}) and  (\ref{eq C6}) shows that $P_{cg}$ satisfies the 
 \emph{Neumann boundary condition}:
 \begin{equation}
\label{eq C8} 
\frac{\partial P_{cg}}{\partial n}=\frac{\partial P_{cg}}{\partial n}^{\hspace{-0,5em}\rm data} \egaldef \left\{ \begin{array}{cl}
   \frac{\sqrt{3}}{3} \,  \nu_{1}^{12,\rm data} \, \displaystyle  \frac{d P_{c}^{12}}{d s_{1}}  
   &    \mbox{(water-oil edge),}    \\
   \frac{\sqrt{3}}{3}\, \big( \nu_{1}^{13,\rm data} \, \displaystyle  \frac{d P_{c}^{12}}{d s_{1}} + \nu_{3}^{13,\rm data}\displaystyle  \frac{d P_{c}^{32}}{d s_{3}}\big)  
    &  \mbox{(water-gas edge),}  \\
         \frac{\sqrt{3}}{3} \,  \nu_{3}^{23,\rm data} \, \displaystyle  \frac{d P_{c}^{32}}{d s_{3}}   
   & \mbox{(gas oil edge),} 
\end{array} \right.
\end{equation}
where $\nu_{1}^{12,\rm data}$,  $ \nu_{3}^{23,\rm data}$, $\nu_{1}^{13,\rm data}$,  $ \nu_{3}^{13,\rm data}$  are defined in (\ref{eq C3-1}) (\ref{eq C3-2}).

Conversely, let $P_{cg}$ satisfy the Dirichlet and Neumann conditions (\ref{eq C4}) (\ref{eq C8}). Then, on the water-oil edge for example, (\ref{eq C4}) and the first line in (\ref{eq C3-1}) implies that 
$\nu_{1}(\CC^{12}(t),p)=\nu_{1}^{12,\rm data}(t,p)$,  and  (\ref{eq C6}) (\ref{eq C7}) and the first line of (\ref{eq C8}) implies that $\nu_{3}(\CC^{12}(t),p)=0$ for $0 \leq t \leq 1$. Hence the fractional flows associated to a global capillary function which satisfies  (\ref{eq C4}) (\ref{eq C8}) honor the data on $\partial \ttt$.

 \subsection{Determination of $d(s,p)$ on $\partial \ttt$.}
 
Now that the global capillary pressure $P_{cg}^{\rm data}$ is known on $\partial \ttt$, honoring the three sets of two-phase relative permeabilities (\ref{eq C1}) on $\partial \ttt$ amounts simply to impose on $d$ the \emph{Dirichlet boundary condition}~:
  
\begin{equation}
\label{eq C13}
  d= d^{\rm data} \quad \mbox{ on } \partial \ttt \ ,
\end{equation}
  where~: 
  \begin{equation}
\label{eq C14} 
d^{\rm data}\egaldef
\left\{ \begin{array}{ll}
    kr_{1}^{12}\,d_{1}(p-P_{cg}^{\rm data} +P_{c}^{12})+kr_{2}^{12}\,d_{2}(p-P_{cg}^{\rm data})  
    												&  \mbox{(water-oil),}  \\ 
    kr_{1}^{13}\,d_{1}(p-P_{cg}^{\rm data} +P_{c}^{12})+kr_{3}^{13}\,d_{3}(p-P_{cg}^{\rm data}+P_{c}^{32})														 & \mbox{(gas-water),} \hspace{-1,5em}  \\
    kr_{3}^{23}\,d_{3}(p-P_{cg}^{\rm data}+P_{c}^{32})+kr_{2}^{23}\,d_{2}(p-P_{cg}^{\rm data}) 
      												&  \mbox{(gas-oil).} 
         \end{array} \right.
\end{equation}

 \subsection{Determination of $P_{cg}$ and $d$ on  the interior of $\ttt$.}

Conditions (\ref{eq C4})(\ref{eq C8})(\ref{eq C14}) are the sole conditions to be satisfied to ensure that the TD-three-phase data derived from $P_{cg}$ and $d$ match the three given sets (\ref{eq C1})(\ref{eq C2}) of two-phase data on $\partial \ttt$ - provided the latter satisfy the TD-compatibility condition (\ref{eq C3-4}).

There is hence a large choice of functions $P_{cg}$ and $d$ to choose from to perform TD-interpolation.

For example, one can use the smoothest functions which satisfy the boundary conditions (\ref{eq C4})(\ref{eq C8})(\ref{eq C14}) by setting:
\begin{equation}
\label{eq C15} \left\{ \begin{array}{ccll}
  \Delta^2 P_{cg}&= &0    & \mbox{in } \ttt   \  , \\
  P_{cg}&=&P_{cg}^{\rm data}    &    \mbox{on }\partial \ttt  \ ,  \\
\displaystyle \frac{\partial P_{cg}}{\partial n}&=&\displaystyle\frac{\partial P_{cg}}{\partial n}^{\hspace{-0,5em}\rm data} 
                                               &  \mbox{on }\partial \ttt \ ,
\end{array}\right.
\end{equation}
and~:
\begin{equation}
\label{eq C16} \left\{ \begin{array}{ccll}
  -\Delta d&=&0    &    \mbox{in } \ttt   \ , \\
  d&=&d^{\rm data}    &     \mbox{on }\partial \ttt  \ .
\end{array} \right.
\end{equation}
These equation can be solved by finite element over the ternary diagram $\ttt$.

One can also try to match some a-priori given three phase target permeability model on the interior of $\ttt$, knowing that the match will not be exact, unless the target relative permeabilities happen to satisfy the TD-compatibility condition (\ref{eq 11-1}). This requires the implementation of an optimization algorithm.

\section[Conclusion]{Conclusion}
\label{sect conclusion}
\setcounter{equation}{0}

\begin{enumerate}
  \item When the three-phase relative permeabilities and capillary pressure satisfy the \emph{TD-condition (\ref{eq 11-1})},  a global capillary pressure function  $P_{cg}(s,p)$ exists such that
  \emph{the classical compressible immiscible three-phase flow equations are fully equivalent to a global pressure formulation}~: the total volumetric flow  $q$ of the three phases and the global pressure $P=P_{2}+P_{cg}(S,P)$ follow the global Darcy law (\ref{eq 10bis}).
  \item  In the global Darcy law (\ref{eq 10bis}), one has always $1-{\partial P_{cg}}/{\partial p}(s,p) >0$   \item If the compressibility of the fluids increases when their wettability decreases (condition (\ref{eq 21})), one has always $1 > 1-{\partial P_{cg}}/{\partial p}(s,p) >0$, and the \emph{stability condition (\ref{eq 8})} is satisfied, which ensures that one can determine $P$ from $P_{2}$ and $S$. Condition (\ref{eq 21})) is always satisfied in water-oil-gas flows where the compressibility of water and oil is neglected.
  \item TD-three-phase relative permeabilities can be obtained by interpolation of three sets of two-phase data on the sides of the ternary diagram, provided the two-phase data satisfy the \emph{TD-compatibility condition (\ref{eq C3-4})}.  TD-interpolation amounts to choose a global capillary function $P_{cg}$ and a total mobility function $d$ which satisfy boundary conditions (\ref{eq C4}) (\ref{eq C8}) (\ref{eq C13}), but can be chosen freely inside the ternary diagram. This interpolation class takes advantage of the lack of information on actual three-phase data inside the ternary diagram to simplify the numerical simulation of the flow.   
\end{enumerate}

\end{document}

%% file: figure_1b_B_W.pstex_t
\begin{picture}(0,0)%
\epsfig{file=figure_1b_B_W.pstex}%
\end{picture}%
\setlength{\unitlength}{4144sp}%
\begingroup\makeatletter\ifx\SetFigFont\undefined%
\gdef\SetFigFont#1#2#3#4#5{%
  \reset@font\fontsize{#1}{#2pt}%
  \fontfamily{#3}\fontseries{#4}\fontshape{#5}%
  \selectfont}%
\fi\endgroup%
\begin{picture}(6628,8127)(-149,-7276)
\put(1915,-1634){\rotatebox{23.0}{\makebox(0,0)[lb]{\smash{{\SetFigFont{20}{24.0}{\rmdefault}{\mddefault}{\updefault}? global pressure ?}}}}}
\put( 91,-3466){\makebox(0,0)[lb]{\smash{{\SetFigFont{14}{16.8}{\rmdefault}{\mddefault}{\updefault}oil saturation}}}}
\put(5956,-3031){\makebox(0,0)[lb]{\smash{{\SetFigFont{14}{16.8}{\rmdefault}{\mddefault}{\updefault}space}}}}
\put(5926,-7231){\makebox(0,0)[lb]{\smash{{\SetFigFont{14}{16.8}{\rmdefault}{\mddefault}{\updefault}space}}}}
\put(-149,-4696){\makebox(0,0)[lb]{\smash{{\SetFigFont{14}{16.8}{\rmdefault}{\mddefault}{\updefault}1}}}}
\put(-149,-7066){\makebox(0,0)[lb]{\smash{{\SetFigFont{14}{16.8}{\rmdefault}{\mddefault}{\updefault}0}}}}
\put(122,-2682){\rotatebox{14.0}{\makebox(0,0)[lb]{\smash{{\SetFigFont{20}{24.0}{\rmdefault}{\mddefault}{\updefault}water pressure}}}}}
\put(4311,203){\rotatebox{9.0}{\makebox(0,0)[lb]{\smash{{\SetFigFont{20}{24.0}{\rmdefault}{\mddefault}{\updefault}oil pressure}}}}}
\put(4666,-5776){\makebox(0,0)[lb]{\smash{{\SetFigFont{20}{24.0}{\rmdefault}{\mddefault}{\updefault}oil}}}}
\put(391,-5656){\makebox(0,0)[lb]{\smash{{\SetFigFont{20}{24.0}{\rmdefault}{\mddefault}{\updefault}water}}}}
\end{picture}%

%% file: figure_3_B_W.pstex_t
\begin{picture}(0,0)%
\epsfig{file=figure_3_B_W.pstex}%
\end{picture}%
\setlength{\unitlength}{4144sp}%
\begingroup\makeatletter\ifx\SetFigFont\undefined%
\gdef\SetFigFont#1#2#3#4#5{%
  \reset@font\fontsize{#1}{#2pt}%
  \fontfamily{#3}\fontseries{#4}\fontshape{#5}%
  \selectfont}%
\fi\endgroup%
\begin{picture}(6670,5510)(-1034,-3636)
\put(2026,1514){\makebox(0,0)[lb]{\smash{{\SetFigFont{29}{34.8}{\rmdefault}{\mddefault}{\updefault}$s_3$}}}}
\put(1801,-2671){\makebox(0,0)[lb]{\smash{{\SetFigFont{29}{34.8}{\rmdefault}{\mddefault}{\updefault}$\mathcal{C}$}}}}
\put(-1034,-3436){\makebox(0,0)[lb]{\smash{{\SetFigFont{29}{34.8}{\rmdefault}{\mddefault}{\updefault}$s_1$}}}}
\put(1481,-561){\makebox(0,0)[lb]{\smash{{\SetFigFont{29}{34.8}{\rmdefault}{\mddefault}{\updefault}${P_{cg}(s)\ ?}$}}}}
\put(3321,-1021){\makebox(0,0)[lb]{\smash{{\SetFigFont{29}{34.8}{\rmdefault}{\mddefault}{\updefault}${s_3}$}}}}
\put(3281,-3501){\makebox(0,0)[lb]{\smash{{\SetFigFont{29}{34.8}{\rmdefault}{\mddefault}{\updefault}${s_1}$}}}}
\put(1881,-1021){\makebox(0,0)[lb]{\smash{{\SetFigFont{29}{34.8}{\rmdefault}{\mddefault}{\updefault}${s}$}}}}
\put(-299,-3446){\makebox(0,0)[lb]{\smash{{\SetFigFont{29}{34.8}{\rmdefault}{\mddefault}{\updefault}water}}}}
\put(4621,-3131){\makebox(0,0)[lb]{\smash{{\SetFigFont{29}{34.8}{\rmdefault}{\mddefault}{\updefault}${0}$}}}}
\put(2321,-1891){\makebox(0,0)[lb]{\smash{{\SetFigFont{29}{34.8}{\rmdefault}{\mddefault}{\updefault}$\ttt$}}}}
\put(1486,779){\makebox(0,0)[lb]{\smash{{\SetFigFont{29}{34.8}{\rmdefault}{\mddefault}{\updefault}gas}}}}
\put(4196,-3481){\makebox(0,0)[lb]{\smash{{\SetFigFont{29}{34.8}{\rmdefault}{\mddefault}{\updefault}oil}}}}
\put(-719,-2791){\makebox(0,0)[lb]{\smash{{\SetFigFont{29}{34.8}{\rmdefault}{\mddefault}{\updefault}${P_{cg}=0}$}}}}
\end{picture}%

%% file: figure_3a_B_W.pstex_t
\begin{picture}(0,0)%
\epsfig{file=figure_3a_B_W.pstex}%
\end{picture}%
\setlength{\unitlength}{4144sp}%
\begingroup\makeatletter\ifx\SetFigFont\undefined%
\gdef\SetFigFont#1#2#3#4#5{%
  \reset@font\fontsize{#1}{#2pt}%
  \fontfamily{#3}\fontseries{#4}\fontshape{#5}%
  \selectfont}%
\fi\endgroup%
\begin{picture}(6670,5510)(-1034,-3636)
\put(1566,-3486){\makebox(0,0)[lb]{\smash{{\SetFigFont{29}{34.8}{\rmdefault}{\mddefault}{\updefault}$\mathcal{C}^1$}}}}
\put(3136,-2136){\makebox(0,0)[lb]{\smash{{\SetFigFont{29}{34.8}{\rmdefault}{\mddefault}{\updefault}$\mathcal{C}^3$}}}}
\put(2026,1514){\makebox(0,0)[lb]{\smash{{\SetFigFont{29}{34.8}{\rmdefault}{\mddefault}{\updefault}$s_3$}}}}
\put(-1034,-3436){\makebox(0,0)[lb]{\smash{{\SetFigFont{29}{34.8}{\rmdefault}{\mddefault}{\updefault}$s_1$}}}}
\put(3281,-3501){\makebox(0,0)[lb]{\smash{{\SetFigFont{29}{34.8}{\rmdefault}{\mddefault}{\updefault}${s_1}$}}}}
\put(1326,-1026){\makebox(0,0)[lb]{\smash{{\SetFigFont{29}{34.8}{\rmdefault}{\mddefault}{\updefault}${s}=(s_1,s_3)$}}}}
\put(1486,779){\makebox(0,0)[lb]{\smash{{\SetFigFont{29}{34.8}{\rmdefault}{\mddefault}{\updefault}gas}}}}
\put(-299,-3446){\makebox(0,0)[lb]{\smash{{\SetFigFont{29}{34.8}{\rmdefault}{\mddefault}{\updefault}water}}}}
\put(4621,-3131){\makebox(0,0)[lb]{\smash{{\SetFigFont{29}{34.8}{\rmdefault}{\mddefault}{\updefault}${0}$}}}}
\put(4196,-3481){\makebox(0,0)[lb]{\smash{{\SetFigFont{29}{34.8}{\rmdefault}{\mddefault}{\updefault}oil}}}}
\put(-719,-2791){\makebox(0,0)[lb]{\smash{{\SetFigFont{29}{34.8}{\rmdefault}{\mddefault}{\updefault}${P_{cg}=0}$}}}}
\put(1736,-2036){\makebox(0,0)[lb]{\smash{{\SetFigFont{29}{34.8}{\rmdefault}{\mddefault}{\updefault}$\ttt$}}}}
\end{picture}%